\documentclass[11pt]{article}
\usepackage{graphicx}
\usepackage{amsmath,amsfonts,amssymb}
\usepackage{citesort}
\usepackage{url}
\usepackage{amsmath}
\usepackage{graphicx}
\usepackage{epsfig}
\usepackage{epstopdf}
\usepackage{color}
\def\tto{\;{\lower 1pt \hbox{$\rightarrow$}}\kern -10pt
\hbox{\raise 2pt \hbox{$\rightarrow$}}\;}

\def\ra{\rangle}
\def\la{\langle}

\def\B{\Bbb B}
\def\h{\hfill\Box}
\def\R{\Bbb R}
\def\N{\Bbb N}
\def\ox{\bar{x}}

\def\dom{\mbox{\rm dom}}

\def\h{\hfill\square}

\def\O{\Omega}
\def\ph{\varphi}

\newcounter{lk}

\begin{document}
\begin{center}
{\bf MINIMIZING DIFFERENCES OF CONVEX FUNCTIONS WITH APPLICATIONS TO FACILITY LOCATION AND CLUSTERING}\\[1ex]
\today\\[2ex]
Nguyen Mau Nam\footnote{Fariborz Maseeh Department of
Mathematics and Statistics, Portland State University, PO Box 751, Portland, OR 97207, United States (email: mau.nam.nguyen@pdx.edu). The research of Nguyen Mau Nam was partially supported by the USA National Science Foundation under grant DMS-1411817.}, Daniel Giles\footnote{Fariborz Maseeh Department of Mathematics and Statistics, Portland State University, PO Box 751, Portland, OR 97207, United States (email: dangiles@pdx.edu). }, R. Blake Rector\footnote{Fariborz Maseeh Department of
Mathematics and Statistics, Portland State University, PO Box 751, Portland, OR 97207, United States (email: r.b.rector@pdx.edu)}.

\end{center}
\small{\bf Abstract.} In this paper we develop algorithms to solve generalized Fermat-Torricelli problems with both positive and negative weights and multifacility location problems involving distances generated by Minkowski gauges. We also introduce a new model of clustering based on squared distances to convex sets. Using the Nesterov smoothing technique and an algorithm for minimizing differences of convex functions called the DCA introduced by Tao and An, we develop effective algorithms for solving these problems.  We demonstrate the algorithms with a variety of numerical examples.\\[1ex]

\newtheorem{Theorem}{Theorem}[section]
\newtheorem{Proposition}[Theorem]{Proposition}
\newtheorem{Remark}[Theorem]{Remark}
\newtheorem{Lemma}[Theorem]{Lemma}
\newtheorem{Corollary}[Theorem]{Corollary}
\newtheorem{Definition}[Theorem]{Definition}
\newtheorem{Example}[Theorem]{Example}
\renewcommand{\theequation}{\thesection.\arabic{equation}}
\normalsize

\section{Introduction}
\setcounter{equation}{0}
The classical Fermat-Torricelli problem asks for a point that minimizes the sum of the Euclidean distances to three points in the plane. This problem was introduced by the French mathematician Pierre De Fermat in the 17th century. In spite of the simplicity of the model, this problem has been a topic for extensive research recently due to both its mathematical beauty and its practical applications in the field of facility location. Several generalized models for the Fermat-Torricelli problem have been introduced and studied in the literature; see \cite{any,b,Tuy92,d,e,HM2015,k,Ku-Ma,Martini,n2,nars,nh,npr,vz} and the references therein.

Given a finite number of target points $a^i\in\R^n$ with the associated weights $c_i\in \R$ for $i=1, \ldots, m$, a generalized model of the Fermat-Torricelli problem seeks to minimize the objective function:
\begin{equation}\label{m1}
f(x):=\sum_{i=1}^m c_i \|x-a^i\|,\qquad \; x\in \R^n.
\end{equation}
Since the weights $c_i$ for $i=1,\ldots,m$ could possibly be negative, the objective function $f$ is not only nondifferentiable but also nonconvex.

A more realistic model asks for a finite number of centroids $x^\ell$ for $\ell=1,\ldots,k$ in $\R^n$  where each $a^i$ is assigned to its nearest centroid. The objective function to be minimized is the weighted sum of the assignment distances:
\begin{equation}\label{m2}
f(x^1, \ldots, x^k):=\sum_{i=1}^m c_i\big (\min_{\ell=1, \ldots, k} \|x^\ell-a^i\|\big),\qquad x^\ell\in \R^n\; \mbox{\rm for }\ell=1, \ldots,k.
\end{equation}
If the weights $c_i$ are nonnegative, \eqref{m1} is a convex function, but \eqref{m2} is nonconvex even if the weights $c_i$ are nonnegative. The problem of minimizing \eqref{m2} reduces to the generalized Fermat-Torricelli problem of minimizing \eqref{m1} in the case where $k=1$. This fundamental problem of multifacility location has a close relationship with clustering problems. Note that the Euclidean distance in objective functions \eqref{m1} and \eqref{m2} can be replaced by generalized distances as necessitated by different applications. Due to the nonconvexity and nondifferentiability of these functions, their minimization needs optimization techniques beyond convexity.

 A recent paper by An, Belghiti, and Tao \cite{abt} used an algorithm called \emph{the DCA} (Difference of Convex Algorithm) to minimize a version of objective function \eqref{m2} that involves the squared Euclidean distances with constant weights $c_i=1$. Their method shows robustness, efficiency, and superiority compared with the well-known $K-$means algorithm when applied to a number of real-world data sets. The DCA was introduced by Tao in 1986, and then extensively developed in the works of An, Tao, and others; see \cite{TA1,TA2} and the references therein. An important feature of the DCA is its simplicity, while still being very effective for many applications compared with other methods. In fact, the DCA is one of the most successful algorithms to deal with nonconvex optimization problems.

In this paper we continue the works of  An, Belghiti, and Tao \cite{abt} by considering the problems of minimizing \eqref{m1} and \eqref{m2} in which the Euclidean distance is replaced by the distance generated by Minkowski gauges. This consideration seems to be more appropriate when viewing these problems as facility location problems. Solving location problems involving Minkowski gauges allows us to unify those generated by arbitrary norms and even more generalized notions of distances; see \cite{HM2015,n2,nars} and the references therein. In addition, our models become nondifferentiable without using squared Euclidean distances as in \cite{abt}. Our approach is based on the Nesterov smoothing technique \cite{n} and the DCA. Based on the DCA, we also propose a method to solve a new model of clustering called \emph{set clustering}. This model involves squared Euclidean distances to convex sets instead of singletons, and hence coincides with the model considered in \cite{abt} when the sets reduce to singletons. Using sets instead of points in  clustering  allows us to classify objects with nonnegligible sizes.

The paper is organized as follows. In Section \ref{sec:IntroDCA}, we give an accessible  presentation of DC programming and the DCA by providing simple proofs for some available results. Section \ref{sec:GenFT} is devoted to developing algorithms to solve generalized weighted Fermat-Torricelli problems involving possibly negative weights and Minkowski gauges. Algorithms for solving multifacility location problems with Minkowski gauges are presented in Section \ref{sec:MultiFac}. In Section \ref{sec:SetCluster} we introduce and develop an algorithm to solve the new model of clustering involving sets.  Finally, we demonstrate our algorithms through a variety of numerical examples in Section \ref{sec:Examples}, and offer some concluding remarks in Section \ref{sec:Conclude}.

\section{An Introduction to the DCA}\label{sec:IntroDCA}
\setcounter{equation}{0}

In this section we provide an easy path to basic results of DC programming and the DCA for the convenience of the reader. Most of the results  in this section can be found in \cite{TA1,TA2}, although our presentation is tailored to the algorithms we present in the following sections.

Consider the problem:
\begin{equation}\label{dcf}
\mbox{\rm minimize}\, f(x):=g(x)-h(x), x\in \R^n,
\end{equation}
where $g\colon\R^n\to (-\infty, \infty]$ and $h\colon \R^n\to \R$ are convex functions. The function $f$ in \eqref{dcf} is called a \emph{DC function} and $g-h$ is called a \emph{DC decomposition} of $f$.

For a convex function $g\colon \R^n\to (-\infty, \infty]$, the \emph{Fenchel conjugate} of $g$ is defined by
\begin{equation*}
g^*(y):=\sup\{\la y, x\ra -g(x)\; |\; x\in \R^n\}.
\end{equation*}
Note that if $g$ is proper, i.e. $\dom(g):=\{x\in \R^n\; |\; g(x)<\infty\}\neq \emptyset$, then $g^*\colon \R^n\to (-\infty, \infty]$ is also a convex function. In addition, if $g$ is lower semicontinuous, then $x\in \partial g^*(y)$ if and only if $y\in \partial g(x)$, where $\partial$ denotes the subdifferential operator in the sense of convex analysis; see, e.g., \cite{HUL,bmn,r}.

The DCA is a simple but effective optimization scheme for minimizing differences of convex functions. Although the algorithm is used for nonconvex optimization problems, the convexity of the functions involved still plays a crucial role. The algorithm is summarized as follows, as applied to (\ref{dcf}).\\[1ex]

{\bf Algorithm 1}.
\begin{center}
\begin{tabular}{| l |}
\hline
{\small INPUT}: $x_1\in \mbox{\rm dom}\,g$, $N\in \N$\\
{\bf for} $k=1, \ldots, N$ {\bf do}\\
\qquad Find $y_k\in \partial h(x_k)$\\
\qquad Find $x_{k+1}\in \partial g^*(y_k)$\\
{\bf end for}\\
{\small OUTPUT}: $x_{N+1}$\\
\hline
\end{tabular}
\end{center}
In what follows, we discuss sufficient conditions for the constructibility of  the sequence $\{x_{k}\}$.

\begin{Proposition}\label{p11} Let $g\colon \R^n\to (-\infty, \infty]$ be a proper lower semicontinuous convex function. Then
\begin{equation*}
\partial g(\R^n):=\bigcup_{x\in \R^n}\partial g(x)=\mbox{\rm dom}\,\partial (g^*):=\{y\in \R^n\; |\; \partial g^*(y)\neq\emptyset\}.
\end{equation*}
\end{Proposition}
{\bf Proof.} Let $x\in \R^n$ and $y\in \partial g(x)$. Then $x\in \partial g^*(y)$ which implies $\partial g^*(y)\neq\emptyset$, and so $y\in\mbox{\rm dom}\,\partial g^*$. The opposite inclusion is just as obvious. $\h$

We say that a function $g\colon \R^n\to (-\infty, \infty]$ is \emph{coercive} if
\begin{equation*}
\lim_{\|x\|\to \infty}\frac{g(x)}{\|x\|}=\infty.
\end{equation*}
We also say that $f$ is \emph{level-bounded} if for any $\alpha\in \R$, the level set $g^{-1}((-\infty, \alpha])$ is bounded.

\begin{Proposition} Let $g\colon \R^n\to (-\infty, \infty]$ be a proper lower semicontinuous convex function. Suppose that $f$ is coercive and level-bounded. Then $\dom(\partial g^*)=\R^n$. In particular, $\dom(g^*)=\R^n$.
\end{Proposition}
{\bf Proof.} It follows from the well-known Br{\o}nsted-Rockafellar theorem that $\partial g(\R^n)$ is dense in $\R^n$; see \cite[Theorem 2.3]{Ruiz}. We first show that  the set $\partial g(\R^n)$ is closed. Fix any sequence $\{v_k\}$ in $\partial g(\R^n)$ that converges to $v$. For each $k\in \N$, choose  $x_k\in \R^n$ such that $v_k\in \partial g(x_k)$. Thus,
\begin{equation}\label{se1}
\la v_k, x-x_k\ra \leq g(x)-g(x_k)\; \mbox{\rm for all }x\in \R^n.
\end{equation}
This implies
\begin{equation*}
g(x_k)-\la v_k, x_k\ra \leq g(x)-\la v_k, x\ra \; \mbox{\rm for all }x\in \R^n.
\end{equation*}
In particular, we can fix $\bar{x}\in \mbox{\rm dom}\,g$ and use the fact that $\{v_k\}$ is bounded to find a constant $\ell_0\in \R$ such that
\begin{equation}\label{se}
g(x_k)-\la v_k, x_k\ra \leq g(\bar{x})-\la v_k, \bar{x}\ra \leq \ell_0\; \mbox{\rm for all }k\in \N.
\end{equation}
Let us now show that $\{x_k\}$ is bounded. By contradiction, assume that this is not the case. Without loss of generality, we can assume that $\lim_{k\to \infty}\|x_k\|=\infty$. By the coercive property of $g$,
\begin{equation*}
\lim_{k\to\infty}\frac{g(x_k)-\la v_k, x_k\ra}{\|x_k\|}=\infty.
\end{equation*}
This is a contradiction to \eqref{se}, so $\{x_k\}$ is bounded. We can assume without loss of generality that $\{x_k\}$ converges to $a\in \R^n$. Then it follows from \eqref{se1} by passing the limit that
\begin{equation*}
\la v, x-a\ra \leq g(x)-g(a)\; \mbox{\rm for all }x\in \R^n.
\end{equation*}
This implies $v\in \partial g(a)\subset \partial g(\R^n)$, and hence $\partial g(\R^n)$ is closed. By Proposition \ref{p11},
\begin{equation*}
\R^n=\partial g(\R^n)=\dom \partial (g^*),
\end{equation*}
which completes the proof. $\h$

Based on the proposition below, we see that in the case where we cannot find $x_k$ or $y_k$ exactly for Algorithm 1, we can find them approximately by solving a convex optimization problem.

\begin{Proposition}\label{characterization} Let $g,h \colon \R^n\to (-\infty, \infty]$ be a proper lower semicontinuous convex function. Then $v\in \partial g^*(y)$ if and only if
\begin{equation}\label{e1}
v\in \mbox{\rm argmin}\,\big\{g(x)-\la y, x\ra \; |\;  x\in \R^n\big\}.
\end{equation}
Moreover, $w\in \partial h(x)$ if and only if
\begin{equation}\label{e2}
w\in \mbox{\rm argmin}\,\big\{h^*(y)-\la y, x\ra \; |\;  y\in \R^n\big\}.
\end{equation}
\end{Proposition}
{\bf Proof.} Suppose that \eqref{e1} is satisfied. Then $0\in \partial \ph(v)$, where $\ph(x):=g(x)-\la y, x\ra$, $x\in \R^n$. It follows that
\begin{equation*}
0\in \partial g(v)-y,
\end{equation*}
and hence $y\in \partial g(v)$ or, equivalently, $v\in \partial g^*(y)$.

Now if we assume that $v\in \partial g^*(y)$, then the proof above gives $0\in \partial \ph(v),$ which justifies \eqref{e1}.

Suppose that \eqref{e2} is satisfied. Then $0\in \partial \psi(w)$, where $\psi(y):=h^*(y)-\la x,y\ra$, $y\in \R^n$. This implies
\begin{equation*}
0\in \partial h^*(w)-x,
\end{equation*}
and hence $x\in \partial h^*(w)$, or, equivalently, $w\in \partial h(x)$. The proof that \eqref{e2} implies $w\in \partial h(x)$ follows as before. $\h$

Based on Proposition \ref{characterization}, we have the another version of the DCA.

{\bf Algorithm 2}.
\begin{center}
\begin{tabular}{| l |}
\hline
{\small INPUT}: $x_1\in \mbox{\rm dom}\,g$, $N\in \N$\\
{\bf for} $k=1, \ldots, N$ {\bf do}\\[1ex]
\qquad Find $y_k\in \partial h(x_k)$ or find $y_k$ approximately by solving the problem:\\[1ex]
\qquad\qquad\qquad\qquad$\mbox{\rm minimize}\; \psi_k(y):=h^*(y)-\la x_k,y\ra,\; y\in \R^n.$\\[1ex]
\qquad Find $x_{k+1}\in \partial g^*(y_k)$ or find $x_{k+1}$ approximately by solving the problem:\\[1ex]
\qquad\qquad\qquad\qquad$\mbox{\rm minimize}\; \phi_k(x):=g(x)-\la x,y_k\ra,\; x \in \R^n.$\\[1ex]
{\bf end for}
{\small OUTPUT}: $x_{N+1}$\\
\hline
\end{tabular}
\end{center}

Let us now discuss the convergence of the DCA.

\begin{Definition} A function $h\colon \R^n\to (-\infty, \infty]$ is called $\gamma$-convex ($\gamma\geq 0$) if there exists $\gamma\geq 0$ such that the function defined by $k(x):=h(x)-\frac{\gamma}{2}\|x\|^2$, $x\in \R^n$, is convex. If there exists $\gamma>0$ such that $h$ is $\gamma-$convex, then $h$ is called strongly convex.
\end{Definition}

\begin{Proposition}\label{gconvex} Let $h\colon \R^n\to (-\infty, \infty]$ be $\gamma$-convex with $\ox\in \mbox{\rm dom}\,h$. Then $v\in \partial h(\ox)$ if and only if
\begin{equation}\label{sc}
\la v, x-\ox\ra +\frac{\gamma}{2}\|x-\ox\|^2\leq h(x)-h(\ox).
\end{equation}
\end{Proposition}
{\bf Proof.} Let $k\colon \R^n\to (-\infty, \infty]$ be the convex function with $k(x)=h(x)-\frac{\gamma}{2}\|x\|^2$. For $v\in \partial h(\ox)$, one has $v\in \partial \ph(\ox)$, where $\ph(x)=k(x)+\frac{\gamma}{2}\|x\|^2$ for $x\in \R^n$. By the subdifferential sum rule,
\begin{equation*}
v\in \partial k(\ox)+\gamma\ox,
\end{equation*}
which implies $v-\gamma\ox\in \partial k(\ox)$. Then
\begin{equation*}
\la v-\gamma\ox, x-\ox\ra \leq k(x)-k(\ox)\; \mbox{\rm for all }x\in \R^n.
\end{equation*}
It follows that
\begin{align*}
\la v, x-\ox\ra &\leq \gamma\la \ox, x\ra -\gamma\la \ox, \ox\ra +h(x)-\frac{\gamma}{2}\|x\|^2-(h(\ox)-\frac{\gamma}{2}\|\ox\|^2)\\
& \leq h(x)-h(\ox)-\frac{\gamma}{2}(\|x\|^2-2\la x, \ox\ra +\|\ox\|^2)\\
&= h(x)-h(\ox)-\frac{\gamma}{2}\|x-\ox\|^2.
\end{align*}
This implies \eqref{sc} and completes the proof. $\h$

\begin{Proposition} Consider the function $f$ defined by \eqref{dcf} and consider the sequence $\{x_k\}$ generated by Algorithm 1. Suppose that $g$ is $\gamma_1$-convex and $h$ is $\gamma_2$-convex. Then
\begin{equation}\label{sc1}
f(x_{k})-f(x_{k+1})\geq \frac{\gamma_1+\gamma_2}{2}\|x_{k+1}-x_k\|^2\; \mbox{\rm for all }k\in \N.
\end{equation}
\end{Proposition}
{\bf Proof.} Since $y_k\in \partial h(x_k)$, by Proposition \ref{gconvex} one has
\begin{equation*}
\la y_k, x-x_k\ra +\frac{\gamma_2}{2}\|x-x_k\|^2\leq h(x)-h(x_k)\; \mbox{\rm for all }x\in \R^n.
\end{equation*}
In particular,
\begin{equation*}
\la y_k, x_{k+1}-x_k\ra +\frac{\gamma_2}{2}\|x_{k+1}-x_k\|^2\leq h(x_{k+1})-h(x_k).
\end{equation*}
In addition, $x_{k+1}\in \partial g^*(y_k)$, and so $y_k\in \partial g(x_{k+1})$, which similarly implies
\begin{equation*}
\la y_k, x_k-x_{k+1}\ra+\frac{\gamma_1}{2}\|x_k-x_{k+1}\|^2 \leq g(x_k)-g(x_{k+1}).
\end{equation*}
Adding these inequalities gives \eqref{sc1}. $\h$

 \begin{Lemma}\label{lm1} Suppose that $h: \R^n\to \R$ is a convex function. If $w_k\in \partial h(x_k)$ and $\{x_k\}$ is a bounded sequence, then $\{w_k\}$ is also bounded.
\end{Lemma}
{\bf Proof.} Fix any point $\ox\in \R^n$. Since $h$ is locally Lipschitz continuous around $\ox$, there exist $\ell>0$ and $\delta>0$ such that
$$|h(x)-h(y)|\leq \ell \|x-y\|\; \mbox{\rm whenever }x,y\in \B(\ox; \delta).$$
This implies that $\|w\|\leq \ell$ whenever $w\in \partial h(u)$ for $u\in \B(\ox; \frac{\delta}{2})$. Indeed,
$$\la w, x-u\ra \leq h(x)-h(u)\; \mbox{\rm for all }x\in \R^n.$$
Choose $\gamma>0$ sufficiently small such that $\B(u; \gamma)\subset \B(\ox; \delta)$. Then
$$\la w, x-u\ra \leq h(x)-h(u)\leq \ell \|x-u\|\; \mbox{\rm whenever }\|x-u\|\leq \gamma.$$
Thus, $\|w\|\leq \ell$.

 For a contradiction, suppose now that $\{w_k\}$ is not bounded. Then we can assume without loss of generality that $\|w_k\|\to \infty$. Since $\{x_k\}$ is bounded, it has a subsequence $\{x_{k_p}\}$ that converges to $x_0\in \R^n$. Let $\ell>0$ be a Lipschitz constant of $f$ around $x_0$. By the observation above,
$$\|w_{k_p}\|\leq \ell\; \mbox{\rm for sufficiently large}\; p.$$
This is a contradiction. $\h$

\begin{Definition} We say that an element $\ox\in \R^n$ is a stationary point of the function $f$ defined by \eqref{dcf} if $\partial g(\ox)\cap \partial h(\ox)\neq \emptyset.$
\end{Definition}

\begin{Theorem} Consider the function $f$ defined by \eqref{dcf} and the sequence $\{x_k\}$ generated by the Algorithm 1. Then $\{f(x_k)\}$ is a decreasing sequence. Suppose further that $f$ is bounded from below, that $g$ is lower semicontinuous, and that $g$ is $\gamma_1$-convex and $h$ is $\gamma_2$-convex with $\gamma_1+\gamma_2>0$. If $\{x_k\}$ is bounded, then every subsequential limit of the sequence $\{x_k\}$ is a stationary point of $f$.
\end{Theorem}
{\bf Proof.} It follows from \eqref{sc1} that $\{f(x_k)\}$ is a decreasing sequence so it converges to real number since $f$ is bounded from below. Then $f(x_k)-f(x_{k+1})\to 0$ as $k\to \infty$, and so using \eqref{sc1} again yields $\|x_{k+1}-x_k\|\to 0$. Suppose that $x_{k_\ell}\to x^*$ as $\ell\to\infty$. By definition,
\begin{equation*}
y_k\in \partial g(x_{k+1})\; \mbox{\rm for all }k\in \N.
\end{equation*}
Since $\{x_k\}$ is bounded, by Lemma \ref{lm1}, $\{y_k\}$ is also a bounded sequence. By extracting a further subsequence, we can assume without loss of generality that $y_{k_\ell}\to y^*$ as $\ell\to \infty$. Since $y_{k_\ell}\in \partial h(x_{k_\ell})$ for all $\ell\in \N$, one has
\begin{equation*}
y^*\in \partial h(x^*).
\end{equation*}
Indeed, by the definition,
\begin{equation}\label{sube}
\la y_{k_\ell}, x-x_{k_\ell}\ra \leq h(x)-h(x_{k_\ell})\; \mbox{\rm for all }x\in \R^n, \ell\in \N.
\end{equation}
Thus,
\begin{equation*}
\la y_{k_\ell}, x^*-x_{k_\ell}\ra \leq h(x^*)-h(x_{k_\ell}).
\end{equation*}
Then $h(x_{k_\ell})\leq \la y_{k_\ell}, x_{k_\ell}-x^*\ra + h(x^*)$, and hence $\limsup h(x_{k_\ell}\ra \leq h(x^*)$. By the lower semicontinuity of $h$, one has that $h(x_{k_\ell})\to h(x^*)$. Letting $\ell\to \infty$ in \eqref{sube} gives $y^*\in \partial h(x^*)$.

Since $\|x_{k+1}-x_k\|\to 0$ and $x_{k_\ell}\to x^*$, one has $x_k\to x^*$. From the relation $y_{k_\ell}\in \partial g(x_{k_{\ell}+1})$, one has $y^*\in \partial g(x^*)$ by a similar argument. Therefore, $x^*$ is a stationary point of $f$. $\h$  

\section{The DCA for a Generalized Fermat-Torricelli Problem}\label{sec:GenFT}
\setcounter{equation}{0}

In this section we develop algorithms for solving the weighted Fermat-Torricelli problem of minimizing \eqref{m1} in which the Euclidean norm is replaced by a Minkowski gauge. Our method is based on the Nesterov smoothing technique and the DCA. This approach allows us to solve generalized versions of the Fermat-Torricelli problem generated by different norms and generalized distances.

Let $F$ be a nonempty closed bounded convex set in $\R^n$ that contains the origin in its interior. Define the \emph{Minkowski gauge} associated with $F$ by
\begin{equation*}
\rho_F(x):=\inf\{t>0\; |\; x\in tF\}.
\end{equation*}
Note that if $F$ is the closed unit ball in $\R^n$, then $\rho_F(x)=\|x\|$.

Given a nonempty bounded set $K$, the \emph{support function} associated with $K$ is given by
\begin{equation*}
\sigma_K(x):=\sup \{\la x, y\ra\; |\; y\in K\}.
\end{equation*}
It follows from the definition of the Minkowski function (see, e.g., \cite[Proposition~2.1]{Ng2006}) that $\rho_F(x)=\sigma_{F^\circ}(x)$, where
\begin{equation*}
F^\circ:=\{y\in \R^n\; |\; \la x, y\ra\leq 1\; \mbox{\rm for all }x\in F\}.
\end{equation*}
Let us present below a direct consequence of the Nesterov smoothing technique given in \cite{n}. In the proposition below, $d(x; \O)$ denotes the Euclidean distance and $P(x; \O)$ denotes the Euclidean projection from a point $x$ to a nonempty closed convex set $\O$ in $\R^n$.
\begin{Proposition}\label{p1} Given any $a\in \R^n$ and $\mu>0$, a Nesterov smoothing approximation of $\ph(x):=\rho_F(x-a)$ has the representation
\begin{equation*}
\ph_\mu(x)=\frac{1}{2\mu}\|x-a\|^2-\frac{\mu}{2}\big[d(\frac{x-a}{\mu}; F^\circ)\big]^2.
\end{equation*}
Moreover, $\nabla \ph_\mu(x)=P(\frac{x-a}{\mu}; F^\circ)$ and
\begin{equation}\label{approx}
\ph_\mu(x)\leq \ph(x)\leq \ph_\mu(x)+\frac{\mu}{2}\|F^\circ\|^2,
\end{equation} where $\|F^\circ\|:=\sup\{\|u\|\; |\; u\in F\}$.
\end{Proposition}
{\bf Proof.} The function $\ph$ can be represented as
\begin{equation*}
\ph(x)=\sigma_{F^\circ}(x-a)=\sup\{\la x-a, u\ra\; |\; u\in F^\circ\}.
\end{equation*}
Using the prox-function $d(x)=\frac{1}{2}\|x\|^2$ in \cite{n}, one obtains a smooth approximation of $\ph$ given by
\begin{align*}
\ph_\mu(x)&:=\sup\{\la x-a, u\ra-\frac{\mu}{2}\|u\|^2\; |\; u\in F^\circ\}\\
&=\sup\{-\frac{\mu}{2}( \|u\|^2-\frac{2}{\mu}\la x-a, u\ra)\; |\; u\in F^\circ\}\\
&=\sup\{-\frac{\mu}{2} \|u-\frac{1}{\mu}(x-a)\|^2+\frac{1}{2\mu}\|x-a\|^2\; |\; u\in F^\circ\}\\
&=\frac{1}{2\mu}\|x-a\|^2-\frac{\mu}{2}\inf \{ \|u-\frac{1}{\mu}(x-a)\|^2\; |\; u\in F^\circ\}\\
&=\frac{1}{2\mu}\|x-a\|^2-\frac{\mu}{2}\big[d(\frac{x-a}{\mu}; F^\circ)\big]^2.
\end{align*}
The formula for computing the gradient of $\ph_\mu$ follows from the well-known gradient formulas for the squared Euclidean norm and the squared distance function generated by a nonempty closed convex set: $\nabla d^2(x; \O)=2[x-P(x; \O)]$; see, e.g., \cite[Exercise 3.2]{bmn}. Estimate \eqref{approx} can be proved directly; see also \cite{n}. The proof is now complete. $\h$

Let $a^i\in \R^n$ for $i=1, \ldots, m$ and let $c_i\neq 0$ for $i=1,\ldots, m$ be real numbers. In the remainder of this section, we study the following generalized version of the Fermat-Torricelli problem:
\begin{equation}\label{FT}
\mbox{\rm minimize}\; f(x):=\sum_{i=1}^m c_i\rho_F(x-a^i),\; x\in \R^n.
\end{equation}
The function $f$ in \eqref{FT} has the following obvious DC decomposition:
\begin{equation*}
f(x)=\sum_{c_i>0}c_i\rho_F(x-a^i)-\sum_{c_i<0}(-c_i)\rho_F(x-a^i).
\end{equation*}
Let $I:=\{i\; |\; c_i>0\}$ and $J:=\{i\; |\; c_i<0\}$ with $\alpha_i=c_i$ if $i\in I$, and $\beta_i=-c_i$ if $i\in J$. Then
\begin{equation}\label{FT1}
f(x)=\sum_{i\in I}\alpha_i \rho_F(x-a^i)-\sum_{j\in J}\beta_j \rho_F(x-a^j).
\end{equation}
Proposition \ref{p2} gives a Nesterov-type approximation for the function $f$.
\begin{Proposition}\label{p2} Consider the function $f$ defined in \eqref{FT1}. Given any $\mu>0$, an approximation of the function $f$ is the following DC function:
\begin{equation*}
f_\mu(x):=g_\mu(x)-h_\mu(x), \; x\in \R^n,
\end{equation*}
where
\begin{align*}
&g_\mu(x):=\sum_{i\in I}\frac{\alpha_i}{2\mu}\|x-a^i\|^2,\\
&h_\mu(x):=\sum_{i\in I}\frac{\mu\alpha_i}{2}\big[d(\frac{x-a^i}{\mu}; F^\circ)\big]^2+\sum_{j\in J}\beta_j \rho_F(x-a^j).
\end{align*}
Moreover, $f_\mu(x)\leq f(x)\leq f_\mu(x)+\frac{\mu \|F^\circ\|^2}{2}\sum_{i\in I}\alpha_i$ for all $x\in \R^n$.
\end{Proposition}
{\bf Proof.} By Proposition \ref{p1},
\begin{align*}
f_\mu(x)&=\sum_{i\in I}\Big[\frac{\alpha_i}{2\mu}\|x-a^i\|^2-\frac{\mu\alpha_i}{2}\big[d(\frac{x-a^i}{\mu}; F^\circ)\big]^2\Big]-\sum_{j\in J}\beta_j \rho_F(x-a_j)\\
&=\sum_{i\in I}\frac{\alpha_i}{2\mu}\|x-a^i\|^2-\Big [\sum_{i\in I}\frac{\mu\alpha_i}{2}\big[d(\frac{x-a^i}{\mu}; F^\circ)\big]^2+\sum_{j\in J}\beta_j \rho_F(x-a^j)\Big].
\end{align*}
The rest of the proof is straightforward. $\h$

\begin{Proposition} Let $\gamma_1:=\sup\{r>0\; |\; B(0; r)\subset F\}$ and $\gamma_2:=\inf\{r>0\; |\; F\subset B(0; r)\}$. Suppose that
\begin{equation*}
\gamma_1\sum_{i\in I}\alpha_i >\gamma_2\sum_{j\in J}\beta_j.
\end{equation*}
Then the function $f$ defined in \eqref{FT1} and its approximation $f_\mu$ defined in Proposition \ref{p2} have absolute minima.
\end{Proposition}
{\bf Proof.} Fix any $r>0$ such that $B(0; r)\subset F$. By the definition, for any $x\in \R^n$,
\begin{equation*}
\rho_F(x)=\inf\{t>0\; |\; t^{-1}x\in F\}\leq \inf\{t>0\; |\; t^{-1}x\in B(0; r)\}=\inf\{t>0\; |\; r^{-1}\|x\|<t\}=r^{-1}\|x\|.
\end{equation*}
This implies $\rho_F(x)\leq \gamma_1^{-1}\|x\|$. Similarly, $\rho_F(x)\geq \gamma_2^{-1}\|x\|$.

Then
\begin{align*}
&\sum_{i\in I}\alpha_i \rho_F(x-a^i)\geq \gamma_2^{-1}\sum_{i\in I}\alpha_i\|x-a^i\|\geq \gamma_2^{-1}\sum_{i\in I}\alpha_i(\|x\|-\|a^i\|),\\
&\sum_{j\in J}\beta_j \rho_F(x-a^j)\leq \gamma_1^{-1}\sum_{j\in J}\beta_j(\|x\|+\|a^j\|).
\end{align*}
It follows that
\begin{equation*}
f(x)\geq \big[(\gamma_2)^{-1}\sum_{i\in I}\alpha_i-(\gamma_1)^{-1}\sum_{j\in J}\beta_j\big]\|x\|-c,
\end{equation*}
where $c:=\gamma_2^{-1}\sum_{i\in I}\alpha_i\|a^i\|+\gamma_1^{-1}\sum_{j\in J}\beta_j\|a^j\|$.

The assumption made guarantees that $\lim_{\|x\|\to \infty}f(x)=\infty$, and so $f$ has an absolute minimum.

By Proposition \ref{p2},
$$f(x)\leq f_\mu(x)+\frac{\mu \|F^\circ\|^2}{2}\sum_{i\in I}\alpha_i.$$
This implies that $\lim_{\|x\|\to \infty}f_\mu(x)=\infty$, and so $f_\mu$ has an absolute minimum as well. $\h$

Define
\begin{equation*}
h^1_\mu(x):=\sum_{i\in I}\frac{\mu\alpha_i}{2}\big[d(\frac{x-a^i}{\mu}; F^\circ)\big]^2, \qquad h^2_\mu(x):=\sum_{j\in J}\beta_j \rho_F(x-a^j) .
\end{equation*}
Then $h_\mu=h^1_\mu+h^2_\mu$ and $h^1_\mu$ is differentiable with
\begin{equation*}
\nabla h^1_\mu(x)=\sum_{i\in I} \alpha_i\big[\frac{x-a^i}{\mu}-P(\frac{x-a^i}{\mu}; F^\circ)\big].
\end{equation*}
\begin{Proposition} Consider the function $g_\mu$ defined in Proposition \ref{p2}. For any $y\in \R^n$, the function
\begin{equation*}
\phi_\mu(x):=g_\mu(x)- \la y, x\ra, x\in \R^n,
\end{equation*}
has a unique minimizer given by
\begin{equation*}
x=\frac{y+\sum_{i\in I}\alpha_i a^i/\mu}{\sum_{i\in I}\alpha_i/\mu}.
\end{equation*}
\end{Proposition}
{\bf Proof.} The gradient of the convex function $\phi_\mu$ is given by
\begin{equation*}
\nabla \phi_\mu(x)=\sum_{i\in I}\frac{\alpha_i}{\mu}(x-a^i)- y.
\end{equation*}
The result then follows by solving $\nabla \phi_\mu(x)=0$. $\h$

Based on the DCA from Algorithm 1, we present the algorithm below to solve the generalized Fermat-Torricelli problem (\ref{FT}):\\[1ex]
{\bf Algorithm 3}.
\begin{center}
\begin{tabular}{| l |}
\hline
{\small INPUTS}: $\mu>0$, $x_1\in \R^n$, $N\in \N$, $F$, $a^1, \ldots, a^m\in \R^n$, $c_1, \ldots, c_m\in \R$.\\[1ex]
{\bf for} $k=1, \ldots, N$ {\bf do}\\[1ex]
\qquad Find $y_k=u_k+v_k$, where\\[1ex]
\qquad\qquad\qquad $u_k:=\sum_{i\in I} \alpha_i\big[\frac{x_k-a^i}{\mu}-P(\frac{x_k-a^i}{\mu}; F^\circ)\big],$\\[1ex]
\qquad\qquad\qquad $v_k\in \sum_{j\in J}\beta_j\partial \rho_F(x_k-a^j).$\\[1ex]
\qquad Find $x_{k+1}=\frac{y_k+\sum_{i\in I}\alpha_i a^i/\mu}{\sum_{i\in I}\alpha_i/\mu}.$\\[1ex]
{\small OUTPUT}: $x_{N+1}$.\\
\hline
\end{tabular}
\end{center}

\begin{Remark}{\rm It is not hard to see that
\begin{equation*}
\partial \rho_F(x)=\left\{
        \begin{array}{ll}
            F^\circ & \quad \mbox{\rm if } x= 0, \\
            \{u\in \R^n\; |\; \sigma_F(u)=1, \la u, x\ra =\rho_F(x)\} & \quad \mbox{\rm if } x \neq 0
        \end{array}
    \right.
    \end{equation*}
    In particular, if $\rho_F(x)=\|x\|$, then
    \begin{equation*}
\partial \rho_F(x)=\left\{
        \begin{array}{ll}
            \B & \quad \mbox{\rm if } x= 0 \\
            \big\{\frac{x}{\|x\|}\big\} & \quad \mbox{\rm if } x \neq 0
        \end{array}
    \right.
    \end{equation*}
    }

\end{Remark}

Let us introduce another algorithm to solve the problem. This algorithm is obtained by using the Nesterov smoothing method for all functions involved in the problem. The proof of the next proposition follows directly from Proposition \ref{p1} as in the proof of Proposition \ref{p2}.

\begin{Proposition}\label{p3} Consider the function $f$ defined in \eqref{FT1}. Given any $\mu>0$, a smooth approximation of the function $f$ is the following DC function:
\begin{equation*}
f_\mu(x):=g_\mu(x)-h_\mu(x), \; x\in \R^n,
\end{equation*}
where
\begin{align*}
&g_\mu(x):=\sum_{i\in I}\frac{\alpha_i}{2\mu}\|x-a^i\|^2,\\
&h_\mu(x):=\sum_{j\in J}\frac{\beta_j}{2\mu}\|x-a^j\|^2-\sum_{j\in J}\frac{\mu\beta_j}{2}\big[d(\frac{x-a^j}{\mu}; F^\circ)\big]^2+\sum_{i\in I}\frac{\mu\alpha_i}{2}\big[d(\frac{x-a^i}{\mu}; F^\circ)\big]^2.
\end{align*}
Moreover, $$f_\mu(x)-\frac{\mu \|F^\circ\|^2}{2}\sum_{i\in I}\beta_i\leq f(x)\leq f_\mu(x)+\frac{\mu \|F^\circ\|^2}{2}\sum_{i\in I}\alpha_i$$ for all $x\in \R^n$.
\end{Proposition}

Note that both functions $g_\mu$ and $h_\mu$ in Proposition \ref{p3} are smooth with the gradients given by
\begin{align*}
 \nabla g_\mu(x)&=\sum_{i\in I}\frac{\alpha_i}{\mu}(x-a^i)\\
 \nabla h_\mu(x)&=\sum_{j\in J}\frac{\beta_j}{\mu}(x-a^j)-\sum_{j\in J}\beta_j\big[\frac{x-a^j}{\mu}-P(\frac{x-a_j}{\mu}; F^\circ)\big]+\sum_{i\in I}\alpha_i\big[\frac{x-a^i}{\mu}-P(\frac{x-a^i}{\mu}; F^\circ)\big]\\
&=\sum_{j\in J}\beta_j\big[P(\frac{x-a^j}{\mu}; F^\circ)\big]+\sum_{i\in I}\alpha_i\big[\frac{x-a^i}{\mu}-P(\frac{x-a^i}{\mu}; F^\circ)\big].
\end{align*}
Based on the DCA in Algorithm 1, we obtain another algorithm for solving problem \eqref{FT}.\\[1ex]
{\bf Algorithm 4}.
\begin{center}
\begin{tabular}{| l |}
\hline
{\small INPUTS}: $\mu>0$, $x_1\in \R^n$, $N\in \N$, $F$, $a^1, \ldots, a^m\in \R^n$, $c_1, \ldots, c_m\in \R$.\\[1ex]
{\bf for} $k=1, \ldots, N$ {\bf do}\\[1ex]
\qquad Find $y_k=u_k+v_k$, where\\[1ex]
\qquad\qquad\qquad $u_k:=\sum_{i\in I}\alpha_i\big[\frac{x_k-a^i}{\mu}-P(\frac{x_k-a^i}{\mu}; F^\circ)\big].$\\[1ex]
\qquad\qquad\qquad $v_k:=\sum_{j\in J}\beta_j\big[P(\frac{x_k-a^j}{\mu}; F^\circ)\big],$\\[1ex]
\qquad Find $x_{k+1}=\frac{y_k+\sum_{i\in I}\alpha_i a^i/\mu}{\sum_{i\in I}\alpha_i/\mu}.$\\[1ex]
{\small OUTPUT}: $x_{N+1}$.\\
\hline
\end{tabular}
\end{center}

\begin{Remark}\label{remark:mu}{\rm When implementing Algorithm 3 and Algorithm 4, instead of using a fixed smoothing parameter $\mu$, we often change $\mu$ during the iteration. The general optimization scheme is
{\small\begin{center}
\begin{tabular}{| l |}
\hline
{\small INITIALIZE}: $x_1\in \R^n$, $\mu_0>0$, $\mu_*>0$, $\sigma\in (0,1)$.\\
Set $k=1$.\\
{\bf Repeat the following}\\
Apply Algorithm 3 (or Algorithm 4) with $\mu=\mu_k$ and starting point $x_k$ \\
\qquad to obtain an approximate solution $x_{k+1}$.\\
\qquad Update $\mu_{k+1}=\sigma \mu_k$.\\
{\bf Until $\mu\leq \mu_*$.}\\
\hline
\end{tabular}
\end{center}}

}
\end{Remark}
%
%
%

\section{Multifacility Location}\label{sec:MultiFac}

In this section we consider a multifacility location problem in which we minimize a general form of the function $f$ defined in \eqref{m2} that involves distances generated by a Minkowski gauge. For simplicity, we consider the case where $c_i=1$ for $i=1,\ldots,m$.

Given $a^i\in \R^n$ for $i=1,\ldots,m$, we need to choose $x^\ell$ for $\ell=1,\ldots, k$ in $\R^n$ as centroids and assign each member $a_i$ to its closest centroid. The objective function to be minimized is the sum of the assignment distances:
\begin{align}\label{mf1}
\mbox{\rm minimize}\, f(x^1, \ldots, x^k)=\sum_{i=1}^m \mbox{\rm min}_{\ell=1,\ldots, k}\, \rho_F(x^\ell-a^i),\qquad x^\ell\in \R^n, \ell=1,\ldots, k.
\end{align}

Let us first discuss the existence of an optimal solution.

\begin{Proposition}\label{existence1} The optimization problem \eqref{mf1} admits  a global optimal solution $(x^1, \ldots, x^k)\in (\R^n)^k$.
\end{Proposition}
{\bf Proof.} We only need to consider the case where $k<m$ because otherwise a global solution can be found by setting $x^\ell=a^\ell$ for $\ell=1,\ldots,m$, and $x^{\ell+1}=\cdots=x^k=a^{m}$. Choose $r>0$ such that
\begin{equation*}
r>\max\{\rho_F(a^i)\; |\; i=1,\ldots, m\}+\max\{\rho_F(a^i-a^j)\; |\; i\neq j\}.
\end{equation*}
Define
\begin{equation*}
\O:=\{(x^1, \ldots, x^k)\in (\R^n)^k\; |\; \rho_F(x^i)\leq r\; \mbox{\rm for all }i=1, \ldots, k\}.
\end{equation*}
Then $\O$ is a compact set. Let us show that
\begin{equation*}
\inf\{f(x^1, \ldots, x^k)\; |\; (x^1, \ldots, x^k)\in \O\}=\inf\{f(x^1, \ldots, x^k)\; |\; (x^1, \ldots, x^k)\in (\R^n)^k\}.
\end{equation*}
Fix any $(x^1, \ldots, x^k)\in (\R^n)^k$. Suppose without loss of generality that $\rho_F(x^i)>r$ for all $i=1, \ldots, p$, where $p\leq k$, and $\rho_F(x^i)\leq r$ for all $i=p+1, \ldots,k$. Since $\rho_F$ is subadditive,
\begin{align*}
\rho_F(x^\ell-a^i)\geq \rho_F(x^\ell)-\rho_F(a^i)> r-\rho_F(a^i)\geq \rho_F(a^\ell-a^i)\; \mbox{\rm for all }\ell=1, \ldots,p, i=1, \ldots, m.
\end{align*}
Therefore,
\begin{align*}
f(x^1, x^2, \ldots, x^k)&=\sum_{i=1}^m \mbox{\rm min}_{\ell=1,\ldots, k}\, \rho_F(x^\ell-a^i)\\
&\geq f(a^1, a^2, \ldots, a^p, x^{p+1}, \ldots, x^k)\geq \inf\{f(x^1, \ldots, x^k)\; |\; (x^1, \ldots, x^k)\in \O\}.
\end{align*}
Thus $\inf\{f(x^1, \ldots, x^k)\; |\; (x^1, \ldots, x^k)\in \O\}\leq \inf\{f(x^1, \ldots, x^k)\; |\; (x^1, \ldots, x^k)\in (\R^n)^k\}$, which completes the proof. $\h$

For our DC decomposition, we start with the following formula:
\begin{equation*}
\mbox{\rm min}_{\ell=1, \ldots, k}\, \rho_F(x^\ell-a^i)=\sum_{\ell=1}^k\rho_F(x^\ell-a^i)-\max_{r=1, \ldots, k}\sum_{\ell=1, \ell\neq r}^k\rho_F(x^\ell-a^i).
\end{equation*}
Then
\begin{equation*}
f(x^1, \ldots, x^k)=\sum_{i=1}^m[\sum_{\ell=1}^k\rho_F(x^\ell-a^i)]-\sum_{i=1}^m\max_{r=1, \ldots, k}\sum_{\ell=1, \ell\neq r}^k\rho_F(x^\ell-a^i)].
\end{equation*}

By Proposition \ref{p1}, the objective function $f$ then has the following approximation:
\begin{equation*}
f_\mu(x^1, \ldots, x^k)=\frac{1}{2\mu}\sum_{i=1}^m\sum_{\ell=1}^k\|x^\ell-a^i\|^2-\big[\frac{\mu}{2}\sum_{i=1}^m\sum_{\ell=1}^k\big[d(\frac{x^\ell-a^i}{\mu}; F^\circ)\big]^2+\sum_{i=1}^m\max_{r=1, \ldots, k}\sum_{\ell=1, \ell\neq r}^k\rho_F(x^\ell-a^i)\big].
\end{equation*}
Thus, $f_\mu(x^1, \ldots, x^k)=g_\mu(x^1, \ldots, x^k)-h_\mu(x^1, \ldots, x^k)$ is a DC decomposition of the function $f_\mu$, where $g_\mu$ and $h_\mu$ are convex functions defined by
\begin{align*}
&g_\mu(x^1, \ldots, x^k):=\frac{1}{2\mu}\sum_{i=1}^m\sum_{\ell=1}^k\|x^\ell-a^i\|^2\; \mbox{\rm and } \\ &h_\mu(x^1, \ldots, x^k):=\frac{\mu}{2}\sum_{i=1}^m\sum_{\ell=1}^k\big[d(\frac{x^\ell-a^i}{\mu}; F^\circ)\big]^2+\sum_{i=1}^m\max_{r=1, \ldots, k}\sum_{\ell=1, \ell\neq r}^k\rho_F(x^\ell-a^i).
\end{align*}

Let $X$ be the $k\times n$-matrix whose rows are $x^1, \ldots, x^k$. We consider the inner product space $\mathcal{M}$ of all $k\times n$ matrices with the inner product of $A, B\in \mathcal{M}$ given by
\begin{equation*}
\la A, B\ra: =\mbox{\rm trace}(AB^T)=\sum_{i=1}^k\sum_{j=1}^n a_{ij}b_{ij}.
\end{equation*}
The norm induced by this inner product is the Frobenius norm.

Then define
\begin{align*}
G_\mu(X):=g_\mu(x^1, \ldots, x^k)&=\frac{1}{2\mu}\sum_{\ell=1}^k\sum_{i=1}^m (\|x^\ell\|^2-2\la x^\ell, a^i\ra +\|a^i\|^2)\\
&=\frac{1}{2\mu}(m\|X\|^2-2\la X, B\ra +k\|A\|^2)\\
&=\frac{m}{2\mu}\|X\|^2-\frac{1}{\mu}\la X, B\ra +\frac{k}{2\mu}\|A\|^2,
\end{align*}
where $A$ is the $m\times n$-matrix whose rows are $a^1, \ldots, a^m$ and $B$ is the $k\times n$-matrix with $a:=\sum_{i=1}^m a^i$ for every row.

Then the function $G_\mu$ is differentiable with gradient given by
$$\nabla G_\mu(X)=\frac{m}{\mu}X-\frac{1}{\mu}B.$$
From the relation $X=\nabla G^*_\mu(Y)$ if and only if $Y=\nabla G_\mu(X)$, one has
$$\nabla G^*_\mu(Y)=\frac{1}{m}(B+ \mu Y).$$

Let us now provide a formula to compute the subdifferential of $H_\mu$ (defined below) at $X$.

Consider the function
\begin{align*}
H^1_\mu(X):&=\frac{\mu}{2}\sum_{i=1}^m\sum_{\ell=1}^k\big[d(\frac{x^\ell-a^i}{\mu}; F^\circ)\big]^2\\
&= \frac{\mu}{2}\Big \{ \big[d(\frac{x^1-a^1}{\mu}; F^\circ)\big]^2+\cdots +[d(\frac{x^1-a^m}{\mu}; F^\circ)\big]^2\Big\}\\
&\cdots\\
& + \frac{\mu}{2}\Big \{ [d(\frac{x^k-a^1}{\mu}; F^\circ)\big]^2+\cdots +[d(\frac{x^k-a^m}{\mu}; F^\circ)\big]^2\Big\}.
\end{align*}
Then the partial derivatives of $H^1_\mu$ are given by
\begin{align*}
&\frac{\partial H^1_\mu}{\partial x^1}(X)=\frac{x^1-a^1}{\mu}-P(\frac{x^1-a^1}{\mu};F^\circ)+\cdots+ \frac{x^1-a^m}{\mu}-P(\frac{x^1-a^m}{\mu};F^\circ)=\sum_{i=1}^m [\frac{x^1-a^i}{\mu}-P(\frac{x^1-a^i}{\mu};F^\circ)],\\
&\vdots\\
&\frac{\partial H^1_\mu}{\partial x^k}(X)=\frac{x^k-a^1}{\mu}-P(\frac{x^k-a^1}{\mu};F^\circ)+\cdots+ \frac{x^k-a^m}{\mu}-P(\frac{x^k-a^m}{\mu};F^\circ)=\sum_{i=1}^m[\frac{x^k-a^i}{\mu}- P(\frac{x^k-a^i}{\mu};F^\circ)].
\end{align*}
The gradient $\nabla H^1_\mu(X)$ is the $k\times n$-matrix whose rows are $\frac{\partial H^1_\mu}{\partial x^1}(X), \ldots, \frac{\partial H^1_\mu}{\partial x^k}(X)$.

Let $H_\mu(X):=h_\mu(x^1, \ldots, x^k)$. Then $H_\mu=H^1_\mu+H^2$, where
\begin{equation*}
H^2(X):=\sum_{i=1}^m\max_{r=1, \ldots, k}\sum_{\ell=1, \ell\neq r}^k\rho_F(x^\ell-a^i).
\end{equation*}
In what follows we provide a formula to find a subgradient of $H^2$ at $X$.

Define the function
\begin{equation*}
F^{i,r}(X):=\sum_{\ell=1, \ell\neq r}^k \rho_F(x^\ell-a^i).
\end{equation*}
Choose  the row vector $v^{i,\ell}\in \partial \rho_F(x^\ell-a^i)$ if $\ell\neq r$ and $v^{i,r}=0$. Then the $k\times n$-matrix formed by the rows $v^{i,r}$ for $i=1, \ldots, k$ is a subgradient of $F^{i,r}$ at $X$.  

Define
\begin{equation*}
F^i(X):=\max_{r=1,\ldots, k} F^{i, r}(X).
\end{equation*}
In order to find a subgradient of $F^i$ at $X$, we first find an index $r\in I_i(X)$, where
\begin{equation*}
I^i(X):=\{r=1, \ldots, k\; |\; F^i(X)=F^{i,r}(X)\}.
\end{equation*}
Then choose $V_i\in \partial F^{i,r}(X)$ and get that $\sum_{i=1}^m V_i$ is a subgradient of the function $H^2$ at $X$.

We have our first algorithm for the multifacility location problem.\\[1ex]
\newpage
{\bf Algorithm 5}.
\begin{center}
\begin{tabular}{| l |}
\hline
{\small INPUTS}: $X_1\in \mathcal{M}$, $N\in \N$, $F$, $a^1, \ldots, a^m\in \R^n$.\\[1ex]
{\bf for} $k=1, \ldots, N$ {\bf do}\\[1ex]
\qquad Find $Y_k=U_k+V_k$, where\\[1ex]
\qquad\qquad\qquad $U_k:=\nabla H^1_\mu(X_k)$,\\[1ex]
\qquad\qquad\qquad $V_k\in \partial H^2(X_k)$.\\[1ex]
\qquad Find $X_{k+1}=\frac{1}{m}(B+\mu Y_k)$.\\[1ex]
{\small OUTPUT}: $X_{N+1}$\\
\hline
\end{tabular}
\end{center}

Let us now present the second algorithm for solving the clustering problem. By Proposition \ref{p1},
 the function $F^{i,r}(X):=\sum_{\ell=1, \ell\neq r}^k\rho_F(x^\ell-a^i)$ has the following smooth approximation of :
 \begin{equation*}
 F_\mu^{i,r}(X)=\sum_{\ell=1, \ell\neq r}^k\big[ \frac{1}{2\mu}\|x^\ell-a^i\|^2-\frac{\mu}{2}\big[d(\frac{x^\ell-a^i}{\mu}; F^\circ)\big]^2\big].
 \end{equation*}
For fixed $r$, define the row vectors $v^{i,\ell}=P(\frac{x^\ell-a^i}{\mu}; F^\circ)$ if $\ell\neq r$ and $v^{i, r}=0$. Then $\nabla F_\mu^{i,r}(X)$ is the $k\times n$ matrix $V_{i,r}$ formed by these rows.

 Now we define the function $F^i_\mu(X):=\max_{r=1, \ldots,k} F_\mu^{i,r}(X)$. This is an approximation of the function
 \begin{equation*}
 F^i(X):=\max_{r=1, \ldots, k}\sum_{\ell=1, \ell\neq r}^k\rho_F(x^\ell-a^i).
 \end{equation*}
 As a result, $H^2_\mu:=\sum_{i=1}^m F^i_\mu$ is an approximation of the function $H^2$.

 Define the active index set
 \begin{equation*}
I^i_\mu(X):=\{r=1, \ldots, k\; |\; F_\mu^i(X)=F_\mu^{i,r}(X)\}.
\end{equation*}
Choose $r\in I_\mu^i(X)$ and calculate $V_i=\nabla F_\mu^{i,r}(X)$. Then $V:=\sum_{i=1}^m V_i$ is a subgradient of the function $H^2_\mu$ at $X$.\\[1ex]
{\bf Algorithm 6}.
\begin{center}
\begin{tabular}{| l |}
\hline
{\small INPUTS}: $X_1\in \mathcal{M}$, $N\in \N$, $F$, $a^1, \ldots, a^m\in \R^n$.\\[1ex]
{\bf for} $k=1, \ldots, N$ {\bf do}\\[1ex]
\qquad Find $Y_k=U_k+V_k$, where\\[1ex]
\qquad\qquad\qquad $U_k:=\nabla H^1_\mu(X_k)$,\\[1ex]
\qquad\qquad\qquad $V_k\in \partial H_\mu^2(X_k)$.\\[1ex]
\qquad Find $X_{k+1}=\frac{1}{m}(B+\mu Y_k)$.\\[1ex]
{\small OUTPUT}: $X_{N+1}$.\\
\hline
\end{tabular}
\end{center}

\begin{Remark}{\rm Similar to the case of Algorithm 3 and Algorithm 4, when implementing Algorithm 5 and Algorithm 6, instead of using a fixed smoothing parameter $\mu$, we often change $\mu$ during the iteration.\\[1ex]

%

}
\end{Remark}

\section{Set Clustering}\label{sec:SetCluster}
\setcounter{equation}{0}

In this section we study the problem of \emph{set clustering}, where the objects being classified are \emph{sets} rather than points. Given a nonempty closed convex set $\O\subset \R^n$, observe that
\begin{align*}
[d(x; \O)]^2&=\inf\{\|x-w\|^2\; |\;  w\in \O\}\\
&=\inf\{\|x\|^2-2\la x, w\ra +\|w\|^2\; |\;  w\in \O\}\\
&=\|x\|^2+\inf\{\|w\|^2-2\la x, w\ra\; |\;  w\in \O\}\\
&=\|x\|^2-\sup\{\la 2 x, w\ra -\|w\|^2\; |\;  w\in \O\}
\end{align*}

\begin{Proposition} Let $\O$ be a nonempty closed  convex set in $\R^n$. Define the function
\begin{equation*}
\ph_{\O}(x):=\sup\{\la 2 x, w\ra -\|w\|^2\; |\;  w\in \O\}=2\sup\{\la  x, w\ra -\frac{1}{2}\|w\|^2\; |\;  w\in \O\}.
\end{equation*}
Then $\ph$ is convex and differentiable with $\nabla \ph_{\O}(x)=2P(x; \O).$
\end{Proposition}
{\bf Proof.} It follows from the representation of $[d(x;\O)]^2$ above that
\begin{equation*}
\ph_{\O}(x)=\|x\|^2-[d(x;\O)]^2.
\end{equation*}
Note that the function $\psi(x):=[d(x;\O)]^2$ is differentiable with $\nabla \psi(x)=2[x-P(x;\O)]$; see, e.g., \cite[Exercise 3.2]{bmn}. Then function $\ph_{\O}$ is differentiable with
\begin{equation*}
\nabla \ph_{\O}(x)=2x-2[x-P(x;\O)]=2P(x;\O),
\end{equation*}
which completes the proof. $\h$

Let $\O^i$ for $i=1,\ldots,m$ be  nonempty closed convex sets in $\R^n$. We need to choose $x^\ell$ for $\ell=1,\ldots, k$ in $\R^n$ as centroids and assign each member $\O^i$ to its closest centroid. The objective function to be minimized is the sum of these distances.

Then we have to solve the optimization problem:
\begin{align}\label{set clustering}
\mbox{\rm minimize}\; f(x^1, \ldots, x^k):=\sum_{i=1}^m \mbox{\rm min}_{\ell=1,\ldots, k}\, [d(x^\ell; \O^i)]^2,\qquad x^\ell\in \R^n, \ell=1,\ldots, k.
\end{align}

\begin{Proposition} Suppose that the convex sets $\O_i$ for $i=1,\ldots,m$ are nonempty closed and bounded. Then \eqref{set clustering} has a global optimal solution.
\end{Proposition}
{\bf Proof.}
Choose $r>0$ such that $\O^i\subset B(0; r)$ for all $i=1, \ldots,m$. Fix $a^i\in \O^i$ for $i=1,\ldots,m$. Define
\begin{equation*}
S:=\{(x^1, \ldots, x^k)\in (\R^n)^k\; |\; \|x^i\|\leq 6r\; \mbox{\rm for }i=1, \ldots, k\}.
\end{equation*}
Let us show that
\begin{equation*}
\inf\{f(x^1, \ldots, x^k)\; |\; (x^1, \ldots, x^k)\in (\R^n)^k\}=\inf\{f(x^1, \ldots, x^k)\; |\; (x^1, \ldots, x^k)\in S\}.
\end{equation*}
Fix any $(x^1, \ldots, x^k)\in (\R^n)^k$. Without loss of generality, suppose that $k<m$ and $\|x^\ell\|>6r$ for $\ell=1,\ldots, p$, and $\|x^{p+1}\|\leq 6r$, \ldots, $\|x^k\|\leq 6r$, where $p\leq k$. Let $p^{\ell,i}:=P(x^\ell; \O^i)$. Then for $\ell=1,\ldots,p$, we have
\begin{align*}
[d(x^\ell, \O^i)]^2&=\|x^\ell-p^{\ell, i}\|^2\\
&=\|x^\ell\|^2-2\la x^\ell, p^{\ell, i}\ra +\|p^{\ell, i}\|^2\\
&\geq \|x^\ell\|^2-2\|x^\ell\| \,\| p^{\ell, i}\|\\
&=\|x^\ell\|(\|x^\ell\|-2\| p^{\ell, i}\|)\geq \|x^\ell\| (6r-2\| p^{\ell, i}\|)\geq 4r\|x^\ell\|\geq 4r^2.
\end{align*}
In addition, for all $\ell=1,...,m$, we have
\begin{equation*}
[d(a^\ell, \O^i)]^2\leq \|a^\ell-a^i\|^2\leq 4r^2\leq [d(x^\ell, \O^i)]^2.
\end{equation*}
It follows that
\begin{align*}
f(x^1, \ldots, x^k)&=\sum_{i=1}^m \mbox{\rm min}_{\ell=1,\ldots, k}\, [d(x^\ell; \O^i)]^2 \\
&\geq f(a^1, \ldots, a^p, x^{p+1}, x^\ell)\\
&\geq \inf\{f(x^1, \ldots, x^k)\; |\; (x^1, \ldots, x^k)\in S\}.
\end{align*}
The rest of the proof follows from the proof of Proposition \ref{existence1}. $\h$

We use the following formula
\begin{equation*}
\mbox{\rm min}_{\ell=1, \ldots, k}\, \, [d(x^\ell; \O^i)]^2=\sum_{\ell=1}^k\, [d(x^\ell; \O^i)]^2-\max_{r=1, \ldots, k}\sum_{\ell=1, \ell\neq r}^k\, [d(x^\ell; \O^i)]^2.
\end{equation*}
Then
\begin{align*}
f(x^1, \ldots, x^k)&=\sum_{i=1}^m\sum_{\ell=1}^k[d(x^\ell; \O^i)]^2-\big[\sum_{i=1}^m\max_{r=1, \ldots, k}\sum_{\ell=1, \ell\neq r}^k[d(x^\ell; \O^i)]^2\big]\\
&=\sum_{i=1}^m\sum_{\ell=1}^k \|x^\ell\|^2-\Big[\sum_{i=1}^m\sum_{\ell=1}^k\ph_{\O^i}(x^\ell)+\sum_{i=1}^m\max_{r=1, \ldots, k}\sum_{\ell=1, \ell\neq r}^k[d(x^\ell; \O^i)]^2\Big].
\end{align*}
Define
\begin{align*}
&g(x^1, \ldots, x^k):=\sum_{i=1}^m\sum_{\ell=1}^k \|x^\ell\|^2\\
&h(x^1, \ldots, x^k):=\sum_{i=1}^m\sum_{\ell=1}^k\ph_{\O^i}(x^\ell)+\sum_{i=1}^m\max_{r=1, \ldots, k}\sum_{\ell=1, \ell\neq r}^k[d(x^\ell; \O^i)]^2.
\end{align*}
We have the DC decomposition $f=g-h$.

For $X\in \mathcal{M}$, define
\begin{equation*}
G(X):=\sum_{i=1}^m\sum_{\ell=1}^k \|x^\ell\|^2=m\|X\|^2.
\end{equation*}
Thus, $\nabla G^*(X)=\frac{1}{2m}(X)$.

Define
\begin{equation*}
H^1(X):=\sum_{i=1}^m\sum_{\ell=1}^k \ph_{\O^i}(x^\ell).
\end{equation*}

Then
\begin{align*}
&\frac{\partial H^1}{\partial x^1}=2P(x^1; \O^1)+\ldots +2P(x^1; \O^m)\\
&\ldots\\
&\frac{\partial H^1}{\partial x^k}=2P(x^k; \O^1)+\ldots +2P(x^k; \O^m)\\
\end{align*}
Then $\nabla H^1(X)$ is the $k\times n$ matrix whose rows are $\frac{\partial H^1}{\partial x^i}$ for $i=1, \ldots,k$.

Let us now present a formula to compute a subgradient of the function
\begin{equation*}
H^2(X)=\sum_{i=1}^m \max_{r=1, \ldots, k}\sum_{\ell=1, \ell\neq r}^k [d(x^\ell; \O^i)]^2.
\end{equation*}

Define
\begin{equation*}
H_2^i(X):=\max_{r=1, \ldots, k}\sum_{\ell=1, \ell\neq r}^k [d(x^\ell; \O^i)]^2=\max_{r=1, \ldots, k}H_2^{i,r},
\end{equation*}
where
\begin{equation*}
H_2^{i,r}:=\sum_{\ell=1, \ell\neq r}^k [d(x^\ell; \O^i)]^2.
\end{equation*}
Consider the following row vectors
\begin{align*}
&v_{i, \ell}:=2(x^\ell-P(x^\ell; \O^i))\; \mbox{\rm if }\ell\neq r\\
&v_{i,r}:=0.
\end{align*}
Then $\nabla H_2^{i,r}$ is the $k\times n$ matrix whose rows are these vectors.

Define the active index set
\begin{equation*}
I^i(X):=\{r=1, \ldots, k\; |\; H_2^{i, r}(X)=H_2^i(X)\}.
\end{equation*}
Choose $r\in I^i(X)$ and let $V_i:=\nabla H_2^{i, r}(X)$. Then $V:=\sum_{i=1}^m V_i$ is a subgradient of $H^2$ at $X$.\\[1ex]
\newpage
{\bf Algorithm 7}.
\begin{center}
\begin{tabular}{| l |}
\hline
{\small INPUTS}: $X_1\in \mathcal{M}$, $N\in \N$, $\O^1, \ldots, \O^m\in \R^n$\\[1ex]
{\bf for} $k=1, \ldots, N$ {\bf do}\\[1ex]
\qquad Find $Y_k=U_k+V_k$, where\\[1ex]
\qquad\qquad\qquad $U_k:=\nabla H^1(X_k)$,\\[1ex]
\qquad\qquad\qquad $V_k\in \partial H^2(X_k)$.\\[1ex]
\qquad Find $X_{k+1}=\frac{1}{2m}(Y_k).$.\\[1ex]
{\small OUTPUT}: $X_{N+1}$.\\
\hline
\end{tabular}
\end{center}

\section{Numerical Implementation}\label{sec:Examples}
\setcounter{equation}{0}

We demonstrate the above algorithms on several problems. All code is written in \textit{MATLAB} and run on an Intel Core i5 3.00 GHz CPU with 8GB RAM. Unless otherwise stated, we use the closed Euclidean unit ball for the set $F$ associated with the Minkowski gauge. In accordance with Remark \ref{remark:mu}, we use $\mu_*=10^{-6}$, decreasing $\mu$ over 3 implementations, each of which runs until $\sum_{\ell=1}^k d(x_j^\ell,x_{j-1}^\ell) < k\cdot 10^{-6}$, where $k$ is the number of centers and $j$ is the iteration counter.  The starting value $\mu_0$ is specified in each example.  

\begin{figure}[!ht]
\begin{center}
\includegraphics[width=8cm]{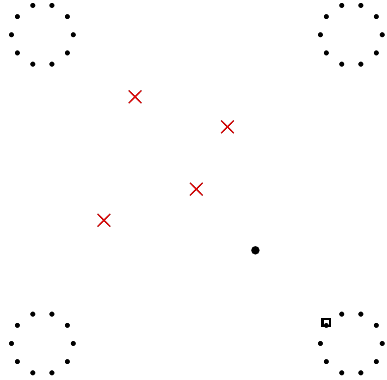}
\caption{ A generalized Fermat-Torricelli problem in $\mathbb{R}^2$. Each $\times$ is a negatively weighted point; the optimal solution is represented by $\bullet$ for the $\ell_2$ norm, and ${\small\pmb{\square}}$ for the $\ell_1$ norm. }
\label{fig:FTex1}
\end{center}
\end{figure}

\textbf{Example 1}\\
In this example we implement Algorithms 3 and 4 to solve a generalized Fermat-Torricelli problem with negative weights, as defined in \eqref{FT}.  We choose $m=44$ points $a_i$ in $\mathbb{R}^2$ as follows. For $i=1,\dots, 40$, we choose distinct $a_i$ from
\[
\{C_h+(\cos(j \pi/5) ,\sin(j\pi/5))\;|\; h=1,\dots,4 \;;\; j=1,\dots,10\}
\]
where each $C_h$ is a distinct element in $\{(\pm 5, \pm 5)\}$ for $h=1,\dots,4$. For $i=41,\dots,44$, we choose distinct $a_i\in\{(0,0),(1,2),(-3,-1),(-2,3)\}$. The weights are assigned $c_i=1$ if $1\leq i\leq 40$ and $c_i=-2$ if $41\leq i \leq 44$.  For the smoothing parameter, we use an initial $\mu_0=.1$.  Then Algorithms 3 and 4 converge to optimal solutions of $x\approx(1.90,-2.00)$ using $F$ as the closed Euclidean unit ball, and $x\approx(4.19, -4.31)$ with $F$ as the closed $\ell_1$ unit ball (Figure \ref{fig:FTex1}).

\begin{figure}[!ht]
\begin{center}
\includegraphics[width=10cm]{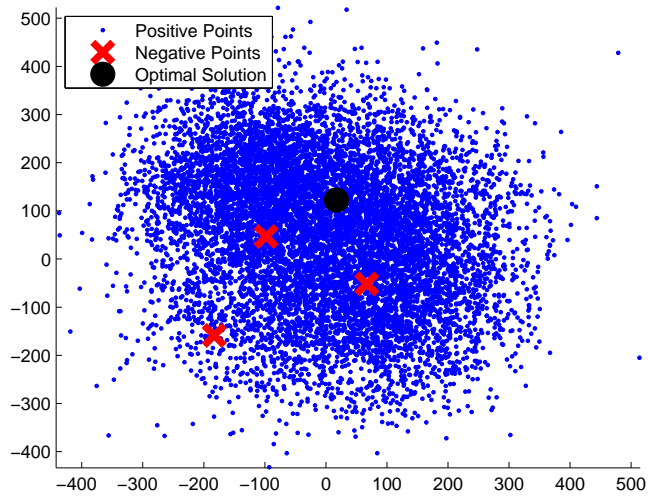}
\caption{ A generalized Fermat-Torricelli problem in $\mathbb{R}^2$. Each negative point has weight of -1000; each positive point has a weight of 1; the optimal solution is represented by $\bullet$ for the $\ell_1$ norm. }
\label{fig:FTex1b_sol}
\end{center}
\end{figure}

\begin{figure}[!ht]
\begin{center}
\includegraphics[width=8cm]{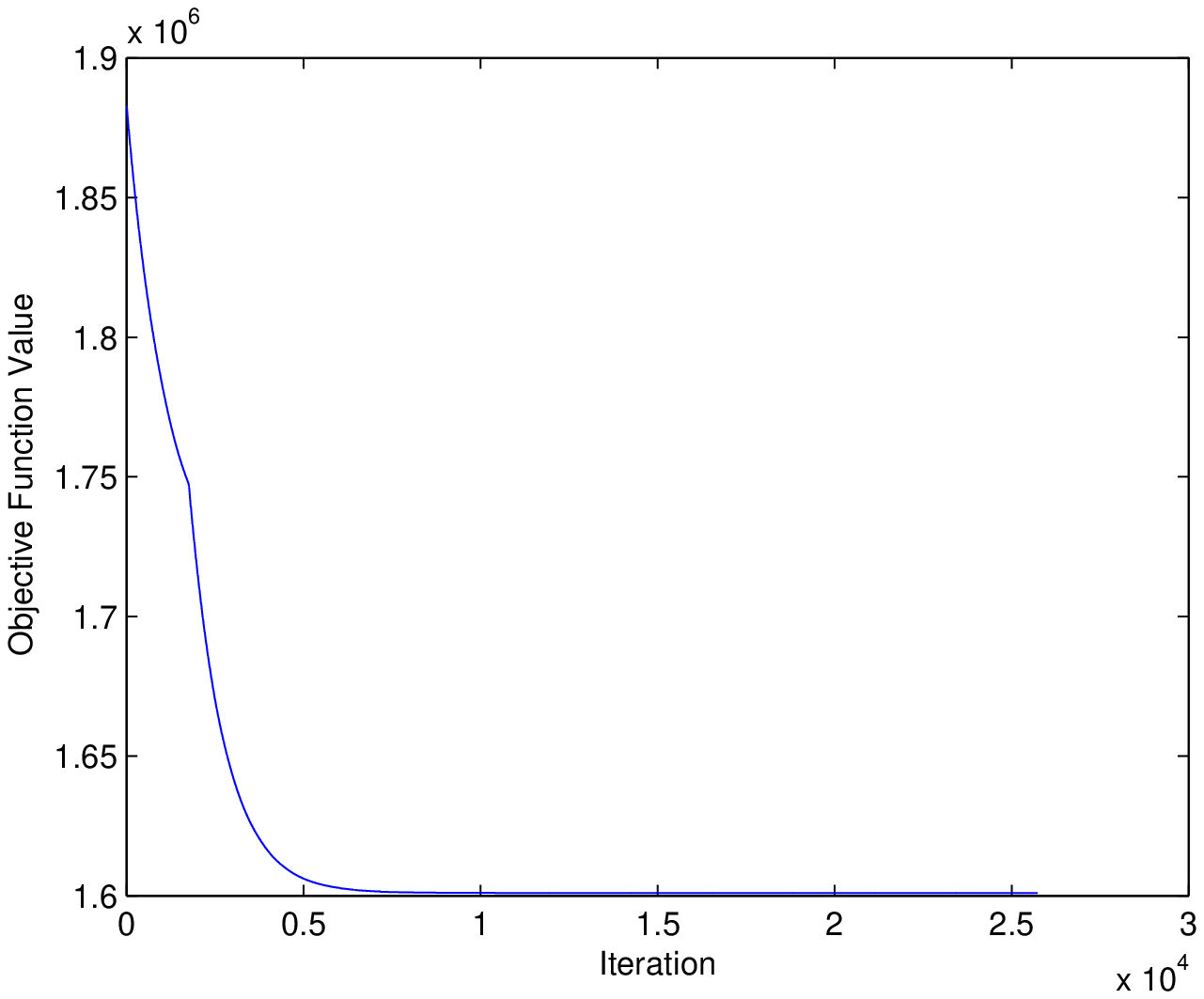}
\caption{ The objective function values for Algorithm 4 for the generalized Fermat-Torricelli problem under the $\ell_1$ norm shown in Figure \ref{fig:FTex1b_sol}.}
\label{fig:FTex1b_convg}
\end{center}
\end{figure}

\textbf{Example 2}\\
In this example we implement Algorithms 4 to solve the generalized Fermat-Torricelli problem under the $\ell_1$ norm with randomly generated points as shown in Figure \ref{fig:FTex1b_sol} .  This synthetic data set has 10,000 points with weight $c_i=1$ and three points with weight $c_i=-1000$.  For the smoothing parameter, we use an initial $\mu_0=.1$.  Then, both Algorithm 4 converges to an optimal solution of $x\approx(17.29,122.46)$.  The convergence rate is shown in Figure \ref{fig:FTex1b_convg} .

\begin{figure}[!ht]
\begin{center}
\includegraphics[width=14cm]{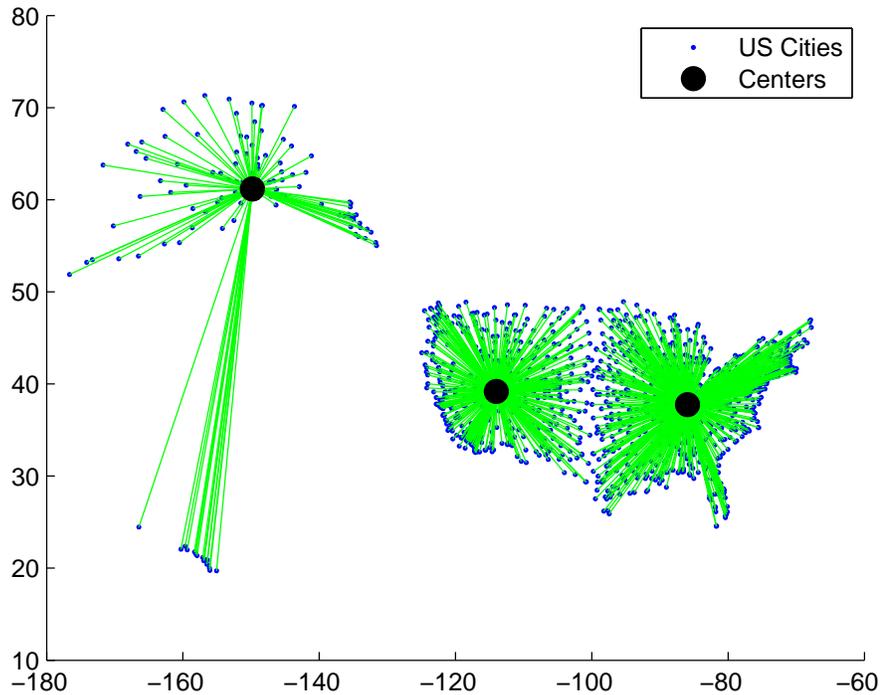}
\caption{ The solution to the multifacility location problem with three centers and Euclidean distance to 1217 US Cities.  A line connects each city with its closest center.}
\label{fig:USCities}
\end{center}
\end{figure}

\textbf{Example 3}\\
We implement Algorithm 5 to solve multifacility location problems given by function (\ref{mf1}). We use the following six real data sets\footnote{Available at \textit{https://archive.ics.uci.edu/ml/datasets.html}}: \textit{WINE} contains 178 instances of $k=3$ wine cultivars in $\mathbb{R}^{13}$. The classical \textit{IRIS} data set contains 150 observations in $\mathbb{R}^4$, describing $k=3$ varieties of Iris flower.  The \textit{PIMA} data set contains 768 observations, each with $8$ features describing the medical history of adults of Pima American-Indian heritage. \textit{IONOSPHERE} contains data on $351$ radar observations in $\mathbb{R}^{34}$ of free electrons in the ionosphere. \textit{USCity}\footnote{http:/www.realestate3d.com/gps/uslatlongdegmin.htm} contains the latitude and longitude of $1217$ US cities; we use $k=3$ centroids (Figure \ref{fig:USCities}).

Reported values are as follows: $m$ is the number of points in the data set; $n$ is the dimension; $k$ is the number of centers; $\mu_0$ is the starting value for the smoothing parameter $\mu$, as discussed in \ref{remark:mu} (in each case, $\sigma$ is chosen so that $\mu$ decreases to $\mu_*$ in three iterations); $Iter$ is the number of iterations until convergence; $CPU$ is the computation time in seconds; $Objval$ is the final value of the true objective function $(1.2)$, not the smoothed version $f_\mu$.  Implementations of Algorithm 6 produced nearly identical results on each example and thus are not reported.

\begin{table}
\begin{center}
\begin{tabular}{l|lllllll|}\cline{2-8}
																		     	& $m$&$n$&$k$ &$\mu_0$& Iter & CPU   & Objval  \\ \hline
\multicolumn{1}{|c}{	\textit{WINE}}  		& 178& 13&  3 & 10    & 690  & 1.86  & $1.62922\cdot 10^4$ \\ \hline
\multicolumn{1}{|c}{	\textit{IRIS}}      & 150& 4 &  3 & 0.1   & 314  & 0.66  & $96.6565$   \\ \hline
\multicolumn{1}{|c}{	\textit{PIMA}}      & 768& 8 &  2 & 10    & 267  & 2.22  & $4.75611\cdot 10^4$  \\ \hline
\multicolumn{1}{|c}{	\textit{IONOSPHERE}}& 351& 34&  2 & 0.1   & 391  & 1.68  & $7.93712\cdot 10^2$  \\ \hline
\multicolumn{1}{|c}{	\textit{USCity}}		&1217& 2 &  3 & 1     & 940  & 16.0 & $1.14211\cdot 10^4$  \\ \hline
	\end{tabular}\\
\caption{Results for Example 2, the performance of Algorithm 5 on real data sets.}
\label{tab:alg5}
\end{center}
\end{table}


\textbf{Example 4}
We now use Algorithm 7 to solve a multifacility location problem involving distances to \emph{sets}, rather than points. We consider the latitude and longitude of the 50 most populous US cities on a plate carr\'{e}e projection.  For demonstration purposes we represent each city with a ball of radius $r=0.1\sqrt{A/\pi}$, where $A$ is the city's reported area in square miles. Coordinates and area for each city were taken from $2014$ United States Census Bureau data\footnote{Available at \textit{https://en.wikipedia.org/wiki/List\_of\_United\_States\_cities\_by\_population}}. Then Algorithm 7 is implemented to minimize function (5.1) with $k=5$ centroids. An optimal solution is given below, and shown in Figure \ref{fig:50cities}.
 \[
	X=\begin{bmatrix}
	 36.2350^\circ N &  77.7130^\circ W \\
	41.1278   ^\circ N&  86.1934^\circ W  \\
	34.2681^\circ N &  95.3486^\circ W  \\
	35.1042^\circ N  &	108.1652 ^\circ W \\
38.2494^\circ N &   120.1098^\circ W   \end{bmatrix}
	\]
\begin{figure}[!ht]
\begin{center}
\includegraphics[scale=1,clip,trim= 0 3cm 0 3cm]{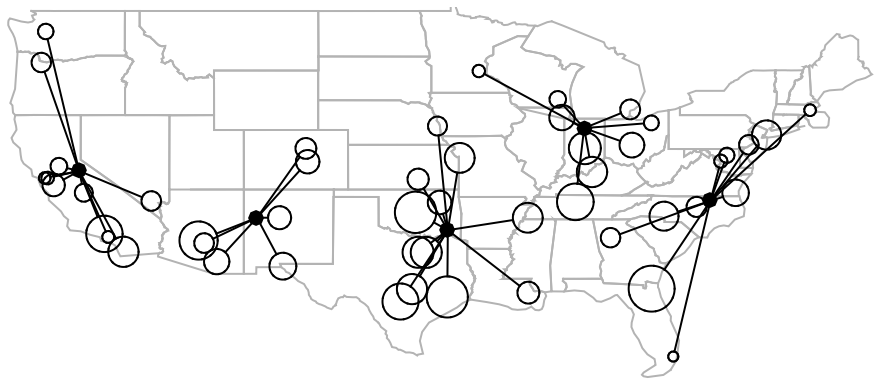}
\caption{The fifty most populous US cities, approximated by a ball proportional to their area. Each is to assigned the closest of five centroids ($\bullet$), which are the optimal facilities. See Example 4.}
\label{fig:50cities}
\end{center}
\end{figure}

\section{Concluding Remarks}\label{sec:Conclude}
\setcounter{equation}{0}

Based on the DCA and the Nesterov smoothing technique, we develop algorithms to solve a number of continuous optimization problems of facility location. Our development continues the works in \cite{abt}. Although unconstrained optimization problems are considered, an easy technique using the \emph{indicator function} and the Euclidean projection would solve the constrained versions of the problems. Another important question is the convergence rate of the algorithms, which can be addressed using recent progress in applying the Kurdyka - Lojasiewicz inequality.


\begin{thebibliography}{99}
 \bibitem{abt} Le Thi Hoai An, M.T. Belghiti, P.D. Tao, A new efficient algorithm based on DC programming and DCA for clustering, J. Glob. Optim.,   27 (2007), 503--608.

\bibitem{amt} L.T.H. An, L.H.  Minh, P.D. Tao, New and efficient DCA based algorithms for minimum sum-of-squares clustering, Pattern Recognition,   47 (2014), 388--401.



\bibitem{any} N.T. An, N.M. Nam, N.D. Yen , A d.c. algorithm via convex analysis approach for solving a location problem involving sets
Journal of Convex Analysis, in press.


\bibitem{b} J. Brimberg, The Fermat Weber location problem revisited, Math. Program. 71 (1995), 71--76.




\bibitem{Tuy92} P.-C. Chen, P. Hansen, B. Jaumard, H. Tuy,
Weber's problem with attraction and repulsion, J. Regional Sci. 32 (1992), 467--486.



\bibitem{d} Z. Drezner, On the convergence of the generalized Weiszfeld algorithm, Ann. Oper. Res.
167 (2009), 327--336.



 \bibitem{e} U. Eckhardt, Weber's problem and Weiszfeld's algorithm in general spaces,  Math. Program. 18 (1980), 186--196.

\bibitem{HM2015} T. Jahn, Y.S. Kupitz, H.  Martini, C. Richter, Minsum location
extended to gauges and to convex sets, J. Optim. Theory Appl. 166 (2015), 711--746.


\bibitem{Ng2006} He, Y.,  Ng, K. F.: Subdifferentials of a minimum time function in Banach spaces, J. Math. Anal.
Appl. 321(2006), 896--910.

\bibitem{HUL} J. B. Hiriart-Urruty and C. Lemar\'{e}chal, Funndamental of Convex Analysis, Springer-Verlag, 2001.

\bibitem{k} H. W. Kuhn, A note on Fermat-Torricelli problem, Math. Program. 4 (1973), 98--107.

  \bibitem{Ku-Ma} Y.S. Kupitz and  H. Martini,  Geometric aspects of the generalized Fermat-Torricelli problem, Bolyai Soc. Math. Stud. 6 (1997), 55--127.




\bibitem{Martini} H. Martini, K.J. Swanepoel, G. Weiss, The Fermat-Torricelli problem in
normed planes and spaces, J. Optim. Theory Appl. 115 (2002), 283--314.

\bibitem{bmn} B. S. Mordukhovich and N. M. Nam, An Easy Path to Convex Analysis and Applications, Morgan \& Claypool Publishers, 2014.

\bibitem{n2} B.S. Mordukhovich and N.M. Nam, Applications of variational analysis to
a generalized Fermat-Torricelli problem. J. Optim. Theory Appl. 148 (2011), 431--454.

\bibitem{nars} N. M. Nam, N. T. An, R. B. Rector, J. Sun, Nonsmooth algorithms and Nesterov's smoothing technique for generalized Fermat-Torricelli problems. SIAM J. Optim. 24 (2014), 1815--1839.


\bibitem{nh} N.M. Nam, N. Hoang, A generalized Sylvester problem and a generalized Fermat-Torricelli problem, J. Convex Anal. 20 (2013), 669–-87.



\bibitem{npr} S. Nickel, J. Puerto, A.M. Rodriguez-Chia,   An
approach to location models involving sets as existing facilities,
Math. Oper. Res. 28 (2003), 693--715.



\bibitem{n} Yu. Nesterov,  Smooth minimization of non-smooth functions, Math.
Program. 103 (2005), 127--152.



\bibitem{nbook} Yu. Nesterov, Introductory lectures on convex optimization. A basic course. Applied Optimization, 87. Kluwer Academic Publishers, Boston, MA, 2004.

\bibitem{Phelps} R.R. Phelps, Convex functions, monotone operators and differentiability, Lecture Notes in
Math. 1364, 2nd Edition, Springer-Verlag, Berlin, 1993.

\bibitem{Ruiz} M. Ruiz Gal\'an, Convex numerical radius, J. Math. Anal. Appl. 361 (2010), 481--491.

    \bibitem{TA1} P.D. Tao, L.T.H. An, Convex analysis approach
to D.C. programming: Theory, algorithms and applications, Acta
Math. Vietnam. 22 (1997), 289--355.

\bibitem{TA2} P.D. Tao,  L.T.H. An, A d.c. optimization algorithm
for solving the trust-region subproblem, SIAM J. Optim. 8
(1998), 476--505.


\bibitem{OR12} J. Orihuelaa, M. Ruiz Gal\'an, A coercive James's weak compactness theorem and
nonlinear variational problems 75 (2012), 598--611.

\bibitem{r} R. T. Rockafellar, Convex Analysis, Princeton University
Press, Princeton, NJ, 1970.

\bibitem{vz} Y. Vardi, C-H. Zhang,  A modified Weiszfeld algorithm for the Fermat-Weber location problem, Math. Program. 90 (2001), Ser. A, 559--566.


\bibitem{w} E. Weiszfeld, Sur le point pour lequel la somme des distances de $n$ points donn\'es est minimum, T$\hat{\mbox{\rm o}}$hoku Math. J. 43 (1937), 355--386.

\bibitem{wf} E. Weiszfeld and F. Plastria, On the point for which the sum of the distances to n given points is minimum, Ann. Oper. Res. 167 (2009), 7--41.

\bibitem{Wg} C. Witzgall, Optimal location of a central facility: mathematical models and concepts. Technical Report 8388, National Bureau of Standards, 1984.




\end{thebibliography}
\end{document}